\documentclass[10pt,twoside,leqno]{amsart}
\usepackage{amssymb,amsbsy,amsmath,amsfonts,amssymb,amscd,epsfig,
times,graphics,color,xypic}

\definecolor{blue}{cmyk}{1.,1.,0.,0.63}
\definecolor{red}{cmyk}{0.,1.,1.,0.63}
\definecolor{green}{cmyk}{1.,0.,1.,0.63}

\usepackage[latin1]{inputenc}
\newcommand{\C}{\mathbb{C}}

\def\dl{{[\![}}\def\dr{{]\!]}}

\newtheorem{lemma}{Lemma}

\begin{document}

%\Large

\title[Constancy of some formal meromorphic maps]{
Note on constancy of some formal 
\\
meromorphic maps
}

\author{Jo\"el Merker}

\address{
D\'epartement de Math\'ematiques et Applications, UMR 8553
du CNRS, \'Ecole Normale
Sup\'erieure, 45 rue d'Ulm, F-75230 Paris Cedex 05, 
France. \ \
{\it Internet}:
{\tt http://www.cmi.univ-mrs.fr/$\sim$merker/index.html}}

\email{merker@dma.ens.fr}

\date{\number\year-\number\month-\number\day}

\maketitle

Let $M$ be a local real analytic or
formal generic submanifold of $\C^n$ ($n
\geqslant 2$) having positive codimension $d \geqslant 1$ and positive
CR dimension $n - d =: m \geqslant 1$ which passes through the origin
and whose extrinsic complexification $\mathcal{ M} \subset \C^n \times
\C^n$ is represented, in appropriate coordinates $(z, w, \zeta, \xi)
\in \C^m \times \C^d \times \C^m \times \C^d$, by $d$ holomorphic
or formal equations (\cite{ ber1996, m2005, mp2006}):
\[
w_j
=
\overline{\Theta}_j(z,\zeta,\xi)
\ \ \ \ \ \ \ \ 
\text{\rm or equivalently:}
\ \ \ \ \ \ \ \ 
\xi_j
=
\Theta_j(\zeta,z,w)
\ \ \ \ \ \ \ \ \ \ \ \ \
{\scriptstyle{(j\,=\,1\,\dots\,d)}}.
\]
Also, introduce the two collections of $m$ complex vector fields:
\[
\footnotesize
\aligned
\mathcal{L}_k
:=
\frac{\partial}{\partial z_k}+
\sum_{j=1}^d\,
\frac{\partial\overline{\Theta}_j}{\partial z_k}
(z,\zeta,\xi)\,
\frac{\partial}{\partial w_j}
\ \ \ \ \ \ \ 
\text{\normalsize\rm and} 
\ \ \ \ \ \ \ 
\underline{\mathcal{L}}_k
:=
{\partial\over\partial \zeta_k}
+
\sum_{j=1}^d\,
\frac{ \partial\Theta_j}{\partial\zeta_k}
(\zeta, z, w)\,
\frac{\partial}{\partial\xi_j},
\endaligned
\]
obviously tangent to $\mathcal{ M}$ with $\big( \mathcal{ L}_1, \dots,
\mathcal{ L}_m \big)$ spanning the foliation $\mathcal{ M} \cap \{
\tau = {\rm cst.} \}$, and $\big( \underline{ \mathcal{ L}}_1, \dots,
\underline{ \mathcal{ L}}_m \big)$ spanning $\mathcal{ M} \cap \{ t =
{\rm cst.} \}$, where $t = (z, w) \in \C^n$ and $\tau = (\tau, \xi)
\in \C^n$ ({\em see} the figure p.~20 of~\cite{ mp2006}).  Some
authors call $M$ of {\sl finite type in the sense of Bloom-Graham}
when the Lie algebra ${\rm Lie} \big( \mathcal{ L}_{ k_1}, \,
\underline{ \mathcal{ L}}_{ k_2} \big)$ generated by all possible Lie
brackets of the $\mathcal{ L}_{ k_1}$ and of the $\underline{
\mathcal{ L}}_{k_2}$ spans $T \mathcal{ M}$ at the origin.
Others call $M$ {\sl minimal}.

Let $H'$ be a maximally real local real analytic or formal submanifold
of $\C^{ n'}$ ($n ' \geqslant 1$) also passing through the origin,
which, after straightening, can be supposed to be represented by the
$n$ equations $z_i' = \bar z_i'$, $i = 1, \dots, n'$. As the outcome
of~\cite{ jlm2008}, formal
holomorphic or meromorphic maps from $M$ into a maximally real, real
analytic or formal submanifold $H' \subset \C^{ n'}$ are all
uninteresting: they must be constant. 
Segre sets (\cite{ ber1996}) can in
fact be fully avoided to check such 
a simple observation, hopefully. It
is clear that studying $n' = 1$ suffices.  The assertion holds for
formal meromorphic maps.

\def\thelemma{\!}\begin{lemma}
{\rm (\cite{ jlm2008}, Main Theorem)} 
Let $f(t) \in \C \dl t \dr$ and $g (t) \in \C\dl t \dr$ be two not
identically zero complex formal power series satisfying:
\[
0
\equiv
f(t)\,\overline{g}(\tau)
-
g(t)\,\overline{f}(\tau),
\]
for $(t, \tau) \in \mathcal{ M}$. If $M$ is of finite type, then
$\frac{ f (t)}{ g(t)}$ is a nonzero real constant.
\end{lemma}

Here, we say for short that $0 \equiv a ( t, \tau)$ for $(t, \tau) \in
\mathcal{ M}$ when $a \big( t, \zeta, \Theta (\zeta, t) \big) \equiv
0$ in $\C \dl t, \zeta\dr$.  As a corollary, choosing $g (t) \equiv
1$, formal power series that are real on $M$ must be constant.

\proof
The most economical arguments begin by applying the differentiations
$\mathcal{ L}_k$, which yields: $0 \equiv \mathcal{L }_kf \cdot
\overline{ g} - \mathcal{ L}_kg \cdot
\overline{f }$.  Eliminating
$\overline{ f}$ from this equation with the help of $0 \equiv f \cdot
\overline{ g} - g \cdot \overline{ f}$, we get $0 \equiv \big( f \cdot
\mathcal{ L}_k g - g \cdot \mathcal{ L}_k f) \cdot \overline{ g}$,
whence:
\[
0
\equiv
f\cdot\mathcal{L}_{k_1}g
-
g\cdot\mathcal{L}_{k_1}f
\ \ \ \ \ \ \ \ \ \ \
\text{\rm and also}
\ \ \ \ \ \ \ \ \ \ \
0
\equiv
f\cdot\underline{\mathcal{L}}_{k_2}g
-
g\cdot\underline{\mathcal{L}}_{k_2}f,
\]
the second family of equations being trivially satisfied.

Suppose now ${\sf R}$ and ${\sf S}$ are two vector fields on $\mathcal{
M}$ with local holomorphic coefficients (in the coordinates of
$\mathcal{ M}$) such that:
\[
0
\equiv
f\cdot{\sf R}g
-
g\cdot{\sf R}f
\ \ \ \ \ \ \ \ \ \ \
\text{\rm and}
\ \ \ \ \ \ \ \ \ \ \
0
\equiv
f\cdot{\sf S}g
-
g\cdot{\sf S}f.
\]
By eliminating $f$ from these two equations, we get $0 \equiv {\sf S}f
\cdot {\sf R} g - {\sf S} g \cdot {\sf R} f$. Next, we apply ${\sf
S}$ to the first and ${\sf R}$ to the second:
\[
\aligned
0
&
\equiv
{\sf S}f\cdot{\sf R}g
+
f\cdot{\sf SR}g
-
{\sf S}g\cdot{\sf R}f
-
g\cdot{\sf SR}f,
\\
0
&
\equiv
{\sf R}f\cdot{\sf S}g
+
f\cdot{\sf RS}g
-
{\sf R}g\cdot{\sf S}f
-
g\cdot{\sf RS}f,
\endaligned
\]
and we subtract:
\[
\aligned
0
&
\equiv
2({\sf S}f\cdot{\sf R}g
-
{\sf S}g\cdot{\sf R}f)
+
f({\sf SR}g-{\sf RS}g)
-
g({\sf SR}f-{\sf RS}f)
\\
&
\equiv
\ \ \ \ \ \ \ \ \ \ \ \ \ \ \ \
0
\ \ \ \ \ \ \ \ \ \ \ \ \ \ \ \
-
f
\cdot[{\sf R},{\sf S}]g
+
g\cdot[{\sf R},{\sf S}]f.
\endaligned
\]
For instance, we may apply this to ${\sf R} = \mathcal{ L}_{ k_1 }$
and ${\sf S} = \underline{ \mathcal{ L} }_{ k_2}$, and again and
again, to conclude that for every iterated Lie bracket, say $\mathcal{
T}$, between the $\mathcal{ L}$'s and the $\underline{ \mathcal{ L}}$'s,
we have:
\[
0
\equiv
f\cdot\mathcal{T}g
-
g\cdot\mathcal{T}f.
\]
Choosing any formal holomorphic coordinates ${\sf x}_1, \dots, {\sf
x}_{ 2m + d}$ on $\mathcal{ M}$, finite-typeness then implies that for
$i = 1, \dots, 2m+d$:
\[
0
\equiv
f\cdot\partial_{{\sf x}_i} g
-
g\cdot\partial_{{\sf x}_i} f
\ \ \ \ \ \ \ \ \ \ \ \
\text{\rm in}\ \ 
\C\dl 
{\sf x}_1,\dots,{\sf x}_{2m+d}
\dr.
\]
Let us expand $f = F \, {\sf x}^\alpha + \cdots$ and $g = G \, {\sf
x}^\beta + \cdots$ with constants $F \neq 0$ and $G \neq 0$, where the
left out terms denote a sum of higher monomials for, say, the {\sl
graded lexicographic order} $<_{\rm grlex}$. Expanding,
we get $0 \equiv FG \, {\sf x}^{ \alpha + \beta - {\bf 1}_i}\,
(\beta_i - \alpha_i) + \cdots$, whence $\alpha = \beta$. After a real
dilation, $F = G = 1$. Next, suppose $f = {\sf x}^\alpha + F \, {\sf
x}^\beta + \cdots$ and $g = {\sf x}^\alpha + G \, {\sf x}^\gamma +
\cdots$ with $F \neq 0$, $G \neq 0$, with $\alpha <_{\rm grlex} \beta
<_{\rm grlex} \infty$ and $\beta \leqslant_{\rm grlex} \gamma
\leqslant_{\rm grlex} \infty$, with of course $G = 0$ in case $\gamma
= \infty$.  If $\beta <_{\rm grlex} \gamma$, we get: $0 \equiv
F(\alpha_i - \beta_i) {\sf x}^{ \alpha + \beta - {\bf 1}_i} + \cdots$,
which is impossible.  So $\gamma = \beta$ and we get: $0 \equiv {\sf
x}^{ \alpha + \beta - {\bf 1}_i} (\alpha_i - \beta_i) (F-G) + \cdots$,
whence $G = F$. Then the induction runs through, giving
$f = g$.
\endproof

\smallskip\noindent{\bf Suggestion.} This generalizes to 
mappings from covering submanifolds of solutions (\cite{ m2007}).

\vfill
\begin{thebibliography}{1}

\bibitem{ber1996}
{\sc Baouendi}, M.S.; {\sc Ebenfelt}, P.; {\sc Rothschild}, L.P.: 
{\em Algebraicity of holomorphic mappings between real algebraic sets
in $\C^n$}, Acta Math. {\bf 177} (1996), no.~2, 225--273.

\bibitem{jlm2008}
{\sc Juhlin}, R.; {\sc Lamel}, B.; {\sc Meylan}, F.: 
{\em Formal meromorphic functions on manifolds of finite type}, 
{\tt arXiv.0803.2103}

\bibitem{me1998} 
{\sc Merker}, J.:
{\em Vector field construction of Segre sets}, 
{\tt arXiv.org/abs/math.CV/9901010}.

\bibitem{m2005}
{\sc Merker}, J.: 
{\em \'Etude de la r\'egularit\'e analytique de l'application de
r\'eflexion CR formelle}, Annales Fac.  Sci. Toulouse, {\bf XIV}
(2005), no.~2, 215--330.

\bibitem{mp2006}
{\sc Merker, J.}; {\sc Porten, P.}:
{\em Holomorphic extension of CR functions, envelopes of holomorphy
and removable singularities}, International Mathematics Research
Surveys, Volume {\bf 2006}, Article ID 28925, 287 pages.  {\tt
math.CV/0701531}

\bibitem{m2007}
{\sc Merker}, J.: 
{\em Lie symmetries and CR geometry}, 118~pp, Journal of Mathematical
Sciences (N.Y.), to appear (2007).

\end{thebibliography}
\end{document}